\documentclass[12pt]{article}
\usepackage{amsfonts}
\usepackage{amssymb}
\newtheorem{theorem}{Theorem}
\newtheorem{corollary}[theorem]{Corollary}
\newtheorem{remark}[theorem]{Remark}

\newtheorem{example}[theorem]{Example}
\newtheorem{lemma}[theorem]{Lemma}
\newtheorem{proposition}[theorem]{Proposition}
\def\be{\begin{equation} }
\def\ee{\end{equation} }
\def\<{\langle}
\def\>{\rangle}

\def\proof{\noindent{\it Proof: }}
\newcommand{\spa}{\mbox{span}}
\newcommand{\grad}{\mbox{grad}}
\newcommand{\R}{\mathbb{R}} 
\newcommand{\Q}{\mathbb{Q}}
\newcommand{\po}{{\hspace*{-1ex}}{\bf .  }}
\def\qed{\ifhmode\unskip\nobreak\fi\ifmmode\ifinner 
\else\hskip5 pt \fi\fi\hbox{\hskip5 pt \vrule width4 pt   
height6 pt depth1.5 pt \hskip 1pt }} 

\begin{document}

\title{Isometric immersions of warped products}
\author{M. Dajczer and T. Vlachos}
\date{}
\maketitle
%\subjclass[2000]{Primary 53C40; Secondary 53B25.}
%\keywords{Warped product, isometric immersion, second fundamental form, warped
%product immersion}

\begin{abstract}{\small
We provide conditions under which an isometric immersion of a
(warped) product of manifolds into a space form must be a (warped) product of
isometric immersions.}
\end{abstract}

A basic problem in the theory of submanifolds is to provide conditions that imply that
an isometric immersion  of a product manifold must be a product of isometric immersions.
The first contribution to that problem was given by Moore \cite{M} under purely intrinsic 
assumptions. Namely, he showed that an isometric immersion  of a  product of 
Riemannian manifolds  
$$
f\colon M^n=M_1^{m_1}\times\ldots\times M_k^{m_k}\to\R^{n+k}
$$
must be a product $f=f_1\times\cdots\times f_k$ of hypersurfaces $f_i\colon
M_i^{m_i}\to\R^{m_i+1}$ if the codimension equals the number of factors and no
factor has an open subset at which all sectional curvatures vanish. 

Moore's local result trivially  fails if higher codimension is allowed. For
instance, one can compose $f$ as above with any (local) isometric
immersion of $\R^{n+k}$ into $\R^{n+k+\ell}$ for $\ell\ge 1$. 
This suggest that for higher codimension a hypothesis of extrinsic 
nature that forbids compositions may be needed.  In that direction, a natural 
assumption is an upper bound on the size of the $s$-nullities of the immersion;
see Example \ref{ex} below. 
The concept of $s$-nullities was introduced in \cite{CD} and has
been playing an increasing role in the study of rigidity questions of
submanifolds; see \cite{CD},  \cite{DF} and \cite{DT3}.

In this paper, we provide conditions under which an isometric immersion of a
product of manifolds must be a  product of immersions. In fact, we consider a much
more general situation.  We allow the ambient space to have any constant 
sectional curvature and the submanifold can be a warped product of manifolds 
instead of just an ordinary Riemannian product. 

We point out that the very special case of a warped product of only two factors 
in codimension two was  solved in \cite{DT2} in great generality, namely, without 
any assumption of extrinsic nature.

\section{Preliminaries}

On a product $M^n=M_0\times M_1\times...\times M_k$ of connected Riemannian manifolds 
we define a new metric by
$$
\< X,Y \>=\< \pi_{0*}X,\pi_{0*}Y\> +
\sum_{i=1}^{k}(\rho_i\circ \pi_0)^2\,\< \pi_{i*}X,\pi_{i*}Y\>
$$
where $\rho_1,\ldots,\rho_k\in C^\infty(M_0)$ are  positive functions and 
$\pi_i\colon M^n\to M_i$  denotes the canonical projection. We call $M^n$ 
endowed with this metric the \emph{warped product} of $M_0,\ldots,M_k$ with 
warping functions $\rho_1,\ldots,\rho_k$ and denote by 
$$
M^n=M_0\times_{\rho_1}M_1\times...\times_{\rho_k}M_k.
$$ 
For simplicity, from now on we identify the tangent bundles $TM_j$ 
with the corresponding tangent distributions to $M^n$.  
\vspace{1ex}

Let $\Q_c^m$ denote a complete and simply connected space form of sectional 
curvature $c$.  If $c\neq 0$, we always view $\Q_c^m$  as a totally umbilical 
hypersurface of Euclidean space $\R^{m+1}$ with  Riemannian or Lorentzian 
signature according to the sign of $c$. Fix a point $q\in\Q_c^m$ and let
$\Q_c^{m_0},\Q_{c_1}^{m_1},\ldots,\Q_{c_k}^{m_k}$, $m=\sum m_j$,
be submanifolds through $q$  such that the first
one is totally geodesic and all the others are totally umbilical with
 mean curvature vectors $z_1,\ldots,z_k$  at $q$ and $\<z_i,z_j\>=-c$ for $i \neq j$.
The \emph{warped product representation}
$$
\Psi:\Q_c^{m_0}\times_{\sigma_1}\Q_{c_1}^{m_1}\times\dots
\times_{\sigma_k}\Q_{c_k}^{m_k}\to \Q_c^m
$$
of $\Q_c^m$ is the map
$$
\Psi(p_0,p_1,\dots,p_k)=p_0+\sum_{i=1}^k\sigma_i(p_0)(p_i-q)
$$
where the functions  $\sigma_i:\Q_c^{m_0}\to \R_+$ are defined as
$$
\sigma_i(p)=\left\{\begin{array}{l}1+\<p-q,a_i\>\;\;\mbox{if}\;\;c=0,
\vspace*{1.5ex}\\
\<p,a_i\>\;\;\;\;\;\;\;\;\;\;\;\;\;\;\mbox{if}\;\;c\neq 0
\end{array} \right.
$$
and satisfy $\sigma_i(q)=1$ with $a_i=cq-z_i$. \newpage

It was shown by N\"olker \cite{N} that 
any isometry of a warped product with $k+1$ factors onto an open dense 
subset of $\Q_c^m$ arises as the restriction of a warped product 
representation as above. 

Given a warped product representation 
\be\label{wp}
\Psi:\Q_c^{m_0}\times_{\sigma_1}\Q_{c_1}^{m_1} \times\dots
\times_{\sigma_k}\Q_{c_k}^{m_k}\to\Q_c^{n+p}
\ee 
and isometric immersions $f_i\colon M_i\to\Q_{c_i}^{m_i}$, 
$0\leq i\leq  k$ with $c_0=c$, the map
\be\label{f}
f=\Psi\circ (f_0\times\ldots \times f_k)\colon
M^n=M_0\times_{\rho_1}M_1\times\ldots\times_{\rho_k}M_k\to
\Q_c^{n+p}
\ee
is an isometric immersion of the warped product manifold
$M^n$ with warping functions $\rho_i=\sigma_i\circ f_0$.

We call $f$ given by (\ref{f}) a \emph{warped product of isometric immersions}.
It is easy to see that its second fundamental form 
$$
\alpha\colon TM\times TM\to N_fM
$$
is \emph{adapted} to the product structure
of $M^n$. This means that
$$
\alpha(X_i,X_j)=0\;\;\mbox{for all}\;\;X_i\in TM_i,\; X_j\in
TM_j,\;\;i\neq j.
$$
 
The following basic result is due to  N\"olker \cite{N}.

\begin{proposition}\po \label{nolker1}  Let $f\colon 
M^n=M_0\times_{\rho_1}M_1\times\ldots\times_{\rho_k}M_k\to \Q_c^{n+p}$ be an 
isometric immersion with adapted second fundamental form.
Then, there is a warped product representation $\Psi$ of $\Q_c^{n+p}$ and isometric 
immersions $f_i\colon M_i\to\Q_{c_i}^{m_i}$, $0\leq i\leq k$, such that 
$f=\Psi\circ (f_0\times\ldots \times f_k)$ is a warped
product of isometric immersions. 
\end{proposition}

We now consider isometric immersions of  a
Riemannian product, i.e., all warping functions are constant. The following 
fact was used by Moore to prove the result discussed in the introduction. 

\begin{corollary}\po\label{moore}
Let $f\colon M^n=M_1\times\ldots\times M_k\to\R^{n+p}$
be an isometric immersion with adapted second fundamental form.
Then $f=f_1\times\ldots \times f_k$ is a product of isometric immersions 
$f_i\colon M_i\to\R^{n_i}$, $1\leq i\leq k$.  
\end{corollary}

For the case $c\neq 0$, we first observe that by fixing a 
point $\bar{p}\in\Q_c^{m_0}$ in the 
warped product representation (\ref{wp}) we obtain the isometric embedding 
$F\colon\Q_{c_1}^{m_1}\times\dots\times\Q_{c_k}^{m_k}\to\Q_c^{n+p}$ with parallel
second fundamental form and flat normal bundle given by 
\be\label{F}
F(p_1,\dots,p_k)=\Psi(\bar{p},p_1,\dots,p_k).
\ee
The following result is due to Reckziegel \cite{Re}.

\begin{proposition}\po \label{nolker2} 
Let $f\colon M^n=M_1\times\ldots\times M_k\to\Q_c^{n+p}$, $c\neq 0$,
be an isometric immersion with adapted second fundamental form. Then, 
there are isometric immersions $f_i\colon M_i\to\Q_{c_i}^{m_i}$, 
$1\leq i\leq  k$, such that 
$$
f=F\circ (f_1\times\cdots\times f_k)
$$ 
where $F$ is given by (\ref{F}). 
\end{proposition}

\section{The main lemma}

Let $\beta\colon V\times V\to W$ be a symmetric bilinear form
where $V$ and $(W,\<\,,\,\>)$ are real vector spaces of finite dimension 
$n$ and $p$, respectively, equipped with inner products.

The $s$-\emph{nullity} $\nu_s$  of $\beta$ for any integer $1\le s\le p$  is
defined by
$$
\nu_s=\max_{U^s\subset W}\dim\{x\in V :
\beta_{U^s}(x,y)=0\;\mbox{for all}\;y\in V\}.
$$
Here $\beta_{U^s}=\pi_{U^s}\circ\beta$
where $U^s$ is any $s$-dimensional subspace of $W$ and $\pi_{U^s}\colon W\to
U^s$ denotes the orthogonal projection.\medskip

Let $R\colon V\times V\times V\times V\to\R$ 
be the multilinear map with the algebraic properties of the curvature tensor 
defined by
$$
R(x,y,z,w)=\<\beta(x,w),\beta(y,z)\>-\<\beta(x,z),\beta(y,w)\>.
$$
\begin{lemma}\po\label{D} Assume that $2p<n$ and  $\nu_s<n-2s$ for 
all $1\leq s\leq p$. Let $V=V_1\oplus V_2$ be an orthogonal splitting  
such that
$$
R(x,y,z,u)=R(x,y,u,v)=R(x,u,v,w)=0
$$
for any $x,y,z\in V_1$ and $u,v,w\in V_2$. Then,
$$ 
S=\spa\,\{\beta(x,y): x\in V_1\;\mathrm{and}\;y\in V_2\}=0.
$$
\end{lemma}

\proof For $x\in V_1$ we denote by  $B_x\colon V_2\to S$  the linear map 
$$ 
B_x(y)=\beta(x,y)
$$ 
and set 
\be\label{des}
D=\ker B_x\subset V_2.
\ee  
Fix $x\in V_1$ such that $B_x$ has \emph{maximal rank},
i.e.,
\be\label{rank}
\mbox{rank}\, B_x\geq \mbox{rank}\, B_y
\ee
for any $y\in V_1$.  Thus, $\dim D\leq\dim\ker B_y$ for any $y\in V_1$. 

We first argue that  
\be\label{one}
D\subset\ker B_y\;\;\mbox{for any}\;\;y\in V_1.
\ee
 From $R(x,y,v,e_j)=0$ we obtain 
\be\label{zero}
\<B_xe_j,B_yv\>=\< B_xv,B_ye_j\>=0,\;\; 1\le j\le\ell,
\ee
where $v\in D$, $y\in V_1$ and $\{e_1,\ldots,e_\ell\}$ is an 
orthonormal basis of $E$ in the orthogonal splitting $V_2=D\oplus E$.

The rank of $B_{x+ty}$ is at most $\ell$ for any $t\in\R$. Therefore, 
the vectors $B_{x+ty}v=tB_yv,\;B_{x+ty}e_j$, $1\le j\le\ell$, are linearly
dependent. Hence, the Gramm determinant of these vectors is an
identically zero polynomial in~$t$. By (\ref{zero}) the term of 
lowest order is $t^2\|B_yv\|^2G$, where $G$ is the Gramm determinant of the 
linearly independent vectors $B_xe_j, 1\le j\le\ell$.
It follows  that $B_yv=0$ for any $y\in V_1$ and $v\in D$,
and this is (\ref{one}).

Next, we prove that
\be\label{two}
\beta(u,v)-\sum_{i,j=1}^\ell\< u,e_i\>\< v,e_{j}\>\beta(e_i,e_j)\in S^\perp
\ee 
for any $u,v\in V_2$. From (\ref{one}) we obtain
$$
\beta(y,v)=\sum_{j=1}^\ell\<v,e_j\> \beta(y,e_j)
$$  
for any $y\in V_1$. Then  $R(y,u,w,v)=0$ yields
$$
\<\beta(u,v),\beta(y,w)\>=\sum_{j=1}^\ell\<v,e_j\>
\<\beta(y,e_j),\beta(u,w)\>
$$
for any $y\in V_1$ and $u,v,w\in V_2$. In particular,
$$
\<\beta(w,u),\beta(y,e_j)\>=\sum_{i=1}^\ell\< u,e_i\>
\<\beta(y,e_i),\beta(e_j,w)\>.
$$
Hence, 
$$
\< \beta(u,v),\beta(y,w)\>=\sum_{i,j=1}^\ell\<u,e_i\>\<v,e_j\>
\<\beta(y,e_i),\beta(e_j,w)\>.
$$
Using
$$
\<\beta(y,e_i),\beta(e_j,w)\>=\<\beta(y,w),\beta(e_i,e_j)\>
$$
we obtain
$$
\<\beta(u,v)-\sum_{i,j=1}^\ell\<u,e_i\>\<v,e_j\> 
\beta(e_i,e_j),\beta(y,w)\>=0
$$
for any $y\in V_1$ and $u,v,w\in V_2$, and this is (\ref{two}).
\vspace{1ex}
 
We have from (\ref{one}) that 
$\beta(u,y)=0$ if $u\in D$ and $y\in V_1$ and  from~(\ref{two}) 
that $\beta_S(u,v)=0$ if $u\in D$ and $v\in V_2$. Therefore, 
$$
\beta_S(u,e)=0\;\;\mbox{if}\;\; u\in D\;\;\mbox{and}\;\; e\in V.
$$
Suppose that $s=\dim S\neq 0$. Then, choose vectors $x_j\in V_j$, $j=1,2$, 
such that $B_{x_j}\colon V_k\to S$, $j\neq k$, has maximal rank. It follows 
from the above that 
$\beta_S(b,e)=0$ for any $b\in \ker B_{x_1}\oplus \ker B_{x_2}$ and 
$e\in V$. Hence, 
$$
\nu_s\geq \dim\ker B_{x_1} + \dim\ker B_{x_2}\geq n-2s,
$$ 
and this contradicts our assumption.\qed

\section{The results}

In this section, we state and prove the results of this paper. The main tool 
for the proofs is the algebraic lemma given in the preceding section. 
\vspace{1,5ex}

We define the $s$-\emph{nullity} $\nu_s(x)$ of an isometric immersion  
$f\colon M^n\to \tilde M^{n+p}$  at a point $x\in M^ n$ for an integer $1\le
s\le p$ as the $s$-nullity of its second fundamental form at that point.
Notice that $\nu_p(x)$ is the standard index of relative nullity of
$f$ at  $x\in M^n$.\vspace{1ex}

We start with the case of a Riemannian products of manifolds
$$
M^n= M_1^{n_1}\times\ldots\times M_k^{n_k}.
$$ 
\begin{theorem}\po\label{product}
Let $f\colon M^n\to\Q_c^{n+p}$ with $2p<n$ be an isometric immersion such that
$\nu_s<n-2s$ for $1\leq s\leq p$ at any point. 
\begin{itemize}
\item [(i)] If $c=0$, then $f=f_1\times\ldots\times f_k$ is a product of
isometric immersions.
\item [(ii)] If $c\neq 0$, then 
$f=F\circ (f_1\times\cdots\times f_k)$ 
where $F$ is given by (\ref{F}).
\end{itemize}
\end{theorem}

\proof We apply Lemma \ref{D} to the second fundamental form of $f$ at 
any point of $M^n$, and conclude  that it must be  adopted to the product 
structure of the manifold. The proof now follows from  Corollary \ref{moore} 
or Proposition \ref{nolker2} according to $c=0$ or $c\neq 0$.
\vspace{1,5ex}\qed

 In the remaining of the paper we consider the case of isometric immersions
of warped product manifolds 
$$
M^n=M_0^{n_0}\times_{\rho_1}M_1^{n_1}\times\ldots\times_{\rho_k}M_k^{n_k}.
$$

Our next result assumes that the warping functions are pairwise 
linearly independent. This condition should not be seen  as a restriction. 
In fact, if two  warping functions are linearly dependent, we may change 
the metric of one of the factors by an homothety in such a way that both 
factors now have the same  warping function and can be put together in a 
new factor.

Before we give the statements for the warped product case, we recall the 
relations between the  Levi-Civita connections and the curvature tensors 
for a warped product metric (left hand side)  and the Riemannian
product metric: 
\be\label{con}
\nabla_XY=\tilde\nabla_XY
+\sum_{j=1}^k(\<X^j,Y^j\>\eta_j-\<X,\eta_j\>Y^j-\<Y,\eta_j\>X^j)
\ee
and
\begin{eqnarray}\label{curv}
R(X,Y)\!\!\!&=&\!\!\!\tilde{R}(X,Y)-\sum_{i,j=1}^k\<\eta_i,\eta_j\>X^i\wedge Y^j\\
\!\!\!\!\!\!\!\!+\!\!\!\!&&\!\!\!\!\!\!\!\!\sum_{j=1}^k[(\nabla_{X^0}
\eta_j-\<\eta_j,X\>\eta_j)\wedge Y^j\nonumber
+ X^j\wedge (\nabla_{Y^0}\eta_j-\<\eta_j,Y\>\eta_j)]
\end{eqnarray}
where $X^j=(\pi_j)_*X$ and $\eta_j=-\grad\log\rho_j$.

\begin{theorem}\po\label{six}
Let $f\colon M^n\to\Q_c^{n+p}$ with $2p<n$ be an isometric immersion
such that $\nu_s<n-2s$ for $1\leq s\leq p$ at any point. 
Assume that the warping functions are pairwise linearly independent 
everywhere. Then $f$ is a warped  product of isometric immersions.
\end{theorem}

\proof 
We have from (\ref{curv}) that
$$
R(X,Y,Z,U)=R(X,Y,U,V)=R(X,U,V,W)=0
$$
for any $X,Y,Z\in TM_0$ and $U,V,W\in TM_1\oplus\cdots\oplus TM_k$.
Hence Lemma \ref{D} applies to the second fundamental form 
of $f$ at any point of $M^n$.  Thus,
\be\label{first}
\alpha(X,V)=0
\ee if $X\in TM_0$ and
$V\in TM_1\oplus\cdots\oplus TM_k$. 

To conclude that the second fundamental form of the immersion must be adopted 
to the product structure, it remains to  show that
\be\label{uv}
\alpha(U,V)=0\;\;\mbox{if}\;\; U\in TM_i\;\;\mbox {and}\;\; V\in TM_j\;\;\mbox{for}\;\; 
i\neq j\;\;\mbox{and}\;\;i,j\ge 1. 
\ee
To see this, first observe that (\ref{con}) gives 
\be\label{nabla}
\nabla_UV=0\;\;\mbox{and}\;\;  \nabla_XU=\nabla_UX=-\<\eta_i,X\>U
\ee
for any  $X\in TM_0$. The Codazzi equation 
$$
(\nabla^\perp_X\alpha)(U,V)=(\nabla^\perp_U\alpha)(X,V)
$$
of $f$ using (\ref{first}) and (\ref{nabla}) yields
$$
(\nabla^\perp_X\alpha)(U,V) = \alpha(\nabla_XU,V).
$$
Since the left hand side is symmetric in $U$ and $V$, 
it follows using (\ref{nabla}) that
$$
\alpha(U,V)\<\eta_i-\eta_j,X\>=0
$$ 
for any  $X\in TM_0$. Thus (\ref{uv}) holds since  $\eta_i-\eta_j\neq 0$ 
if $i\neq j$ by assumption, and the proof follows from  Proposition \ref{nolker1}. 
\qed

\begin{example}\label{ex}\po
{\em The assumption on the $s$-nullities in Theorem \ref{six} goes
beyond excluding the case of compositions as discussed in the introduction. 
In fact, for submanifolds with two factors in codimension two it was shown 
in \cite{DT2} that there are two families of submanifolds that are not warped 
products of isometric immersions. Although the submanifolds belonging to one 
family are compositions, the ones  in the other family are not.} 
\end{example}

As in Moore's result discussed in the introduction, we next restrict the 
codimension to the number of factors and assume a curvature condition.

\begin{theorem}\po\label{seven} 
Let $f\colon M^n=M_0^{n_0}\times_{\rho_1}M_1^{n_1}
\times\ldots\times_{\rho_k}M_k^{n_k}\to\Q_c^{n+k+1}$  with $2k+2<n$  be an
isometric immersion such that $\nu_s<n-2s$ for $1\leq s\leq k+1$ at any 
point. Assume that no factor has an open subset where
the sectional curvature or the corresponding warping function is constant. 
Then $f$ is a warped product of hypersurfaces. 
\end{theorem}

\proof   
If the  number $r$ of pairwise linearly independent warping functions is $r=k$, 
the result follows from Theorem \ref{six} and the curvature assumption. 
Thus, we may assume $r<k$ and  
let $\rho_{i_1},\dots,\rho_{i_r}$, \mbox{$1\leq i_1<\dots<i_r\leq k$},  be the
pairwise linearly independent warping functions. 
Hence, we may view $M^n$ as a warped product
\be\label{pros}
M^n=M_0^{\ell_0}\times_{\rho_{i_1}}\hat{M}_1^{\ell_1}\times\cdots
\times_{\rho_{i_r}}\hat{M}_r^{\ell_r}
\ee
where the factors $\hat{M}_j^{\ell_j}$ are the Riemannian products
$$
\hat{M}_j^{\ell_j}=\Pi_{i\in I_j}M_i,\;1\leq j\leq r,
$$
and $I_j$ denotes the set of all indices $1\leq i\leq k$ that correspond to
factors with the same associated warping function after homoteties, if necessary. 

We apply Theorem \ref{six} to $M^n$ with the warped product 
structure (\ref{pros}). It follows that $f$ is a warped product of 
isometric immersions 
$$
f=\Psi\circ (\hat{f_0}\times\hat{f}_1\times\cdots\times\hat{f}_r)
$$ 
with respect to a warped product representation
$$
\Psi\colon\Q_{c}^{m_0}\times_{\sigma_1}\Q_{c_1}^{m_1}\times\cdots
\times_{\sigma_r}\Q_{c_r}^{m_r}\to\Q_c^{n+k+1},
$$ 
where
$\hat{f_0}=f_0$ and
$\hat{f}_j\colon\hat{M}_j^{\ell_j}\to\Q_{c_j}^{m_j},\;\; 1\leq j\leq r.$

We show that $\hat{f}_j$, $1\leq j\leq r$,  satisfies at any 
point of $\hat{M}_j^{\ell_j}$ that 
\be\label{claim}
\nu^{\hat{f}_j}_s<\ell_j-2s
\ee
for $1\leq s\leq\mathrm{cod}(\hat{f}_j)=m_j-\ell_j$. 
We have from Lemma $12$ in \cite{N} that 
\be\label{al}
\alpha=\Psi_*\alpha^{\hat{f}_0\times\cdots\times\hat{f}_r},
\ee
where 
$$
\alpha^{\hat{f}_0 \times\cdots\times \hat{f}_r}(X,Y)=
(\gamma(X,Y),\alpha^{\hat{f}_1}(\hat{X}_1,\hat{Y}_1),\ldots,
\alpha^{\hat{f}_r}(\hat{X}_r,\hat{Y}_r))
$$
and
$$
\gamma(X,Y)=\alpha^{\hat{f}_0}(X_0,Y_0)-\sum_{j=1}^r\rho_{i_j}
\langle\hat{X}_j,\hat{Y}_j\rangle
(\grad\,\sigma_j-{f_0}_*\,\grad\,\rho_{i_j}).
$$
We argue by contradiction. At some point  of $\hat{M}_{j_0}^{\ell_{j_0}}$, let 
$\nu^{\hat{f}_{j_0}}_{s_0} \geq \ell_{j_0}-2s_0$
for some $1\leq j_0\leq r$  and  $1\leq s_0\leq\mathrm{cod}(\hat{f}_{j_0})$. 
Therefore, there exists a subspace 
$U^{s_0}\subset N_{\hat{f}_{j_0}}\hat{M}_{j_0}$ such that
$$
\dim\,\{Y\in TM_{j_0}:\alpha^{\hat{f}_{j_0}}_{U^{s_0}}(Y,Z)=0\;\;
\mbox{for all}\;\; Z\in TM_{j_0}\}\geq \ell_{j_0}-2s_0.
$$
 From the above, we have
$$
\alpha_{U^{s_0}}^{\hat{f}_0\times\cdots\times\hat{f}_r}(X,Y)
=(0,0,\ldots,\alpha^{\hat{f}_{j_0}}_{U^{s_0}}(\hat{X}_{j_0},\hat{Y}_{j_0}),\ldots,0).
$$
It follows that
$$
\nu_{s_0}(\hat{f}_0\times\cdots \times\hat{f}_r)\geq l_{j_0}-2s_0+\sum_{i\neq j_0}l_i=n-2s_0,
$$
In view of (\ref{al}) this a contradiction and proves (\ref{claim}). 

Assume that  $|I_j|>1$ for some $j$. From (\ref{claim}) and Theorem \ref{product}, we obtain
$$
\hat{f}_j=\left\{\begin{array}{l}
g^j_1\times\cdots\times g^j_{|I_j|}\;\;\;\;\;\;\;\;\;\;\;\;\;\mbox{if}\;\;c_j=0,
\vspace*{1.5ex}\\
F_j\circ(g^j_1\times\cdots\times g^j_{|I_j|})\;\;\;\mbox{if}\;\;c_j\neq 0.
\end{array} \right.
$$
By the curvature assumption $\mathrm{cod}({g}^j_i)\geq 1\;\mbox{for any}\; i\in I_j.$
Therefore, 
$$
\mathrm{cod}(\hat{f}_j)\geq |I_j|\;\; \mbox{if}\;\;c_j=0\;\;\;\mbox{and}\;\;
\mathrm{cod}(\hat{f}_j)>|I_j|\;\;\mbox{if}\;\; c_j\neq 0.
$$
We also have from the curvature assumption that $\mathrm{cod}(\hat{f}_j)\geq 1$
if either $j=0$ or $|I_j|=1$. Hence,
$$
k+1=\sum_{j=0}^r\mathrm{cod}(\hat{f}_j)\geq \sum_{j=0}^r |I_j|=k+1.
$$
Therefore $\mathrm{cod}(\hat{f}_j)=|I_j|$ 
for any $0\leq j\leq r$. In particular, if $|I_j|>1$ then $\hat{f}_j$ is a 
product of Euclidean hypersurfaces. 
We conclude that each factor in the initial product decomposition of $M^n$ must
be a hypersurface.\qed

\begin{remark}\po {\em Notice that the proof of the last result gives 
some additional information. For instance, if the warping functions
are not pairwise linearly independent then we must have $c\leq 0$.
}\end{remark}

The results of this paper are local in nature.  Global results for isometric
immersions of Riemannian products where obtained in \cite{AM} and \cite{BDT}.
In the latter the role of compositions was clarified.

{\renewcommand{\baselinestretch}{1}

\hspace*{-20ex}\begin{tabbing} \indent\= IMPA -- Estrada Dona Castorina, 110
\indent\indent\= Univ. of Ioannina -- Math. Dept. \\
\> 22460-320 -- Rio de Janeiro -- Brazil  \>
45110 Ioannina -- Greece \\
\> E-mail: marcos@impa.br \> E-mail: tvlachos@uoi.gr
\end{tabbing}}

\end{document}